# GAUSSIAN LIMITS FOR GENERALIZED SPACINGS

BY YU. BARYSHNIKOV, MATHEW D. PENROSE AND J. E. YUKICH[1]

*Bell Laboratories, University of Bath and Lehigh University*

Nearest neighbor cells in $R^d, d \in \mathbb{N}$, are used to define coefficients of divergence ($\phi$-divergences) between continuous multivariate samples. For large sample sizes, such distances are shown to be asymptotically normal with a variance depending on the underlying point density. In $d = 1$, this extends classical central limit theory for sum functions of spacings. The general results yield central limit theorems for logarithmic $k$-spacings, information gain, log-likelihood ratios and the number of pairs of sample points within a fixed distance of each other.

**1. Introduction.** Suppose $X_{(i)}, 1 \leq i \leq n$, are the order statistics drawn from an i.i.d. sample with distribution $F$ on $\mathbb{R}$ and let $G$ be a distribution function. Classical spacing functionals on $\mathbb{R}$ (Section 6 of [35]) take the form of an *empirical $\phi$-divergence*

$$\sum_{i=1}^{n-1} \phi(n[G(X_{(i+1)}) - G(X_{(i)})]), \tag{1.1}$$

where $\phi : \mathbb{R}^+ \to \mathbb{R}$ is a specified function and where typically $F$ is unknown. When $F$ and $G$ have densities $f$ and $g$, respectively, the functionals (1.1) represent an empirical version of the *$\phi$-divergence* of $g$ from $f$, namely $\int f(x)\phi(\frac{g(x)}{f(x)})dx$. The $\phi$-divergence functional, introduced by Ali and Silvey [1, 2, 3] and independently by Csiszár [9, 10, 11] is a measure of the discrepancy of $G$ relative to $F$. Empirical $\phi$-divergences are widely used in nonparametric estimation and are well suited for goodness-of-fit tests [7, 8, 15, 22, 34, 36, 42].

This paper has two main goals. The first is to use $k$th nearest neighbor cells to establish high-dimensional analogs of the $\phi$-divergences (1.1). The

Received July 2007; revised February 2008.
[1]Supported by NSA Grant H98230-06-1-0052.

nearest neighbor cells are employed to define the statistical discrepancy of a proposed distribution with density $g$ relative to an observed i.i.d. sample drawn from a distribution with density $f$. We establish a general central limit theorem (CLT) showing that the resulting distance functionals converge to a normal random variable whenever $f$ and $g$ are bounded away from zero and infinity. The limiting variance is given in terms of the $V_{\phi,k}$-divergence and $\Delta_{\phi,k}$-divergence of $g$ from $f$, where $V_{\phi,k}$ and $\Delta_{\phi,k}$ are certain integral transforms of $\phi$.

Our second goal is to use $\phi$-divergences based on $k$th nearest neighbors cells to provide a unifying approach toward proving classical central limit theorems for sum functions of $k$-spacings [7, 8, 12, 13, 18, 40, 42]. This yields asymptotic normality for information gain, log-likelihood ratios and sums of logarithmic spacings whenever the densities of $F$ and $G$ are bounded away from zero and infinity. The methods extend to yield a central limit theorem for the number of pairs of sample points within a fixed distance.

More generally, we consider the natural random measures associated with the empirical $\phi$-divergences, obtained for $d=1$ by putting a point mass at each $X_{(i)}$ of size equal to the $i$th term in (1.1), and analogously for $d>1$. We show that these point measures, when acting on bounded test functions, and when suitably centered and scaled, converge weakly to a Gaussian field.

Our approach uses *stabilization* methods, a tool [5, 31, 32, 33] for establishing general limit theorems for sums of weakly dependent terms in geometric probability. These methods quantify local dependence in ways useful for establishing thermodynamic and Gaussian limits and they also show that locally defined functionals of Poisson points on large bounded sets can be well approximated by globally defined functionals of homogeneous Poisson points on all of $\mathbb{R}^d$. This latter feature conveniently often leads to explicit thermodynamic and variance asymptotics.

Existing general limit results cannot be applied directly to the high-dimensional analogs of (1.1). However, it turns out that the empirical $\phi$-divergences nonetheless involve sums of stabilizing functionals, and one might thus expect that the underlying ideas and methods at the heart of stabilization are applicable and lead to variance asymptotics and Gaussian limits for the high-dimensional analogs of (1.1). This paper shows that this is indeed the case. Further, by adapting the methods of [32] to the present setting, we may prove variance asymptotics and central limit theorems over point sets with a fixed (non-Poisson) number of points.

Our general results are stated in Section 2; applications associated with particular choices of $\phi$ are discussed in Section 3. Because of their generality, our main results have lengthy proofs, which we provide (in Section 4) with some details omitted; for full details, see the extended version of this article [4].

## 2. Main results.

### 2.1. *Preliminaries.*

*Notation.* We use the following notation throughout. If $B$ is a Borel subset of $\mathbb{R}^d$, then $|B|$ denotes its Lebesgue measure. Given $\mathcal{X} \subset \mathbb{R}^d$, $a \geq 0$ and $y \in \mathbb{R}^d$, let $y + a\mathcal{X} := \{y + ax : x \in \mathcal{X}\}$. For $x \in \mathbb{R}^d$, let $|x|$ be its Euclidean modulus and for $r > 0$, let $B_r(x)$ denote the open Euclidean ball $\{y \in \mathbb{R}^d : |y - x| < r\}$. Let $\mathbf{0}$ denote the origin of $\mathbb{R}^d$, and let $\omega_d := |B_1(\mathbf{0})| = \pi^{d/2}/\Gamma((d/2) + 1)$. We use $\log x$ to denote the natural logarithm of $x$.

We let $f$ and $g$ denote two probability density functions on $\mathbb{R}^d$ ($d \in \mathbb{N}$) with common compact support, which we assume is convex and which is denoted by $A$. *We assume once and for all that $f$ and $g$ are bounded and that they are bounded away from zero on $A$.* Abusing notation, we let $F(\cdot)$ [resp. $G(\cdot)$] denote the probability measure on $\mathbb{R}^d$ with density $f$ (resp. $g$), that is, $F(B) := \int_B f(x)\,dx$ and $G(B) := \int_B g(x)\,dx$.

Throughout $X_1, X_2, \ldots$ denotes a sequence of independent random $d$-vectors with common density $f$. Let $\mathcal{X}_n := \{X_1, \ldots, X_n\}$. Also, given $\lambda > 0$, let $\mathcal{P}_\lambda$ be a Poisson point process in $A$ with intensity function $\lambda f : A \to \mathbb{R}^+$. For all $a > 0$, let $\mathcal{H}_a$ denote a homogeneous Poisson point process on $\mathbb{R}^d$ with intensity $a$. We write $\mathcal{H}$ for $\mathcal{H}_1$.

Given a Borel subset $E \subset \mathbb{R}^d$, let $\mathcal{B}(E)$ denote the class of bounded Borel-measurable real-valued functions on $E$. Given $h \in \mathcal{B}(\mathbb{R}^d)$, we write $\|h\|_\infty$ for $\sup_{x \in \mathbb{R}^d}(|h(x)|)$ and given also $\mu$ a Borel measure on $\mathbb{R}^d$, we let $\langle h, \mu \rangle$ denote the integral of $h$ with respect to $\mu$.

We shall consider $\phi$-divergences and related quantities for a general class $\mathcal{F}$ of functions $\phi$, which we now describe. Let $\mathbb{R}^+ := (0, \infty)$. Given a continuous function $\phi : \mathbb{R}^+ \to \mathbb{R}$, define the function $\phi^* : \mathbb{R}^+ \to [0, \infty)$ by

$$(2.1) \qquad \phi^*(t) := \begin{cases} \sup\{|\phi(u)| : t \leq u \leq 1\}, & \text{if } 0 < t \leq 1, \\ \sup\{|\phi(s)| : 1 \leq s \leq t\}, & \text{if } t \geq 1. \end{cases}$$

In other words, $\phi^*$ is the minimal function on $\mathbb{R}^+$ with the properties that (i) $-\phi^*(\cdot)$ is unimodal with a maximum at 1, and (ii) $\phi^*(\cdot)$ dominates $|\phi(\cdot)|$ pointwise.

Let $\mathcal{F}$ be the class of continuous functions $\phi : \mathbb{R}^+ \to \mathbb{R}$ such that the restriction to $(0, 1)$ of the function $\phi^*$ defined by (2.1) is square-integrable on $(0, 1)$, and such that $\log(\max(\phi(t), 1)) = o(t)$ as $t \to \infty$. Let $\mathcal{F}_0$ be the class of functions in $\mathcal{F}$ which are bounded on $(0, 1]$.

Let $\Gamma_1$ denote a gamma$(1, 1)$ random variable, that is, let $\Gamma_1$ be exponentially distributed with mean one. Letting $\Gamma_{1,i}$, $i \geq 1$, be independent copies of $\Gamma_1$, we put $\Gamma_k := \sum_{i=1}^k \Gamma_{1,i}$, a gamma random variable with parameters $k$ and 1. For $\sigma^2 > 0$, let $N(0, \sigma^2)$ denote a normal random variable with mean zero and variance $\sigma^2$. Given random variables $X, Y$ we write $X \prec Y$ (or $Y \succ X$) if $Y$ dominates $X$ stochastically, that is, if $P[X \leq t] \geq P[Y \leq t]$ for all $t \in \mathbb{R}$.

2.2. *High-dimensional $\phi$-divergence based on $k$-nearest neighbor cells.* Let $\mathcal{K}$ be an open convex cone in $\mathbb{R}^d$ (a cone is a set that is invariant under dilations). For all $r > 0$, let $B_r^{\mathcal{K}}(x) := x + (\mathcal{K} \cap B_r(\mathbf{0}))$. Recall that the aspect ratio of a subset $E$ of $\mathbb{R}^d$ is the ratio of the radius of the smallest ball containing $E$ and the radius of the largest ball contained in $E$. For $d \geq 2$, we assume that $\mathcal{K}$ is "regular" with respect to $A$, that is, $\mathcal{K}$ is chosen such that the sets $B_r^{\mathcal{K}}(x) \cap A$ have bounded aspect ratio uniformly over $x \in A, r > 0$. When $\mathcal{K} = \mathbb{R}^d$, this condition is trivially satisfied. If $A$ is the unit cube, then $\mathcal{K}$ may be either a tilted orthant or a right circular cone not tangent to any coordinate subspace.

Given the cone $\mathcal{K}$, $x \in \mathbb{R}^d$, a finite set $\mathcal{X} \subset \mathbb{R}^d$, and $k \in \mathbb{N}$, put

$$(2.2) \quad C_k(x, \mathcal{X}) := C_k^{\mathcal{K}}(x, \mathcal{X}) := \bigcup_{t > 0 \,:\, \mathrm{card}(B_t^{\mathcal{K}}(x) \cap \mathcal{X} \setminus \{x\}) < k} B_t^{\mathcal{K}}(x).$$

Here, $\mathrm{card}(\mathcal{Y})$ denotes the cardinality of the finite set $\mathcal{Y}$. If $\mathrm{card}((x + \mathcal{K}) \cap \mathcal{X} \setminus \{x\}) \geq k$, then $C_k^{\mathcal{K}}(x, \mathcal{X})$ is the largest set of the form $B_t^{\mathcal{K}}(x)$ containing fewer than $k$ points of $\mathcal{X} \setminus \{x\}$; otherwise, $C_k^{\mathcal{K}}(x, \mathcal{X})$ is the whole "wedge" $x + \mathcal{K}$. When $\mathcal{K} = \mathbb{R}^d$, $C_k^{\mathcal{K}}(x, \mathcal{X})$ is a ball whose radius is the distance between $x$ and its $k$th nearest neighbor in $\mathcal{X} \setminus x$.

For each $n \geq 2$ and $X_i, 1 \leq i \leq n$, we use the *directed nearest neighbor cells* $C_k^{\mathcal{K}}(X_i, \mathcal{X}_n)$ to define high-dimensional spacing functionals analogous to the classical one-dimensional functionals (1.1). Define for $1 \leq i \leq n$ the transformed $k$th nearest neighbor spacings

$$D_{i,n,k}^g := G(C_k^{\mathcal{K}}(X_i, \mathcal{X}_n)).$$

Given $\phi \in \mathcal{F}$, define the random point measure $\nu_{n,\phi,k}^g$, with total measure $N_{n,\phi,k}^g$, as follows:

$$(2.3) \quad \nu_{n,\phi,k}^g := \sum_{i=1}^n \phi(n D_{i,n,k}^g) \delta_{X_i}; \qquad N_{n,\phi,k}^g := \sum_{i=1}^n \phi(n D_{i,n,k}^g).$$

Here, $\delta_x$ denotes the unit point mass at $x$. Let $\bar{\nu}_{n,\phi,k}^g := \nu_{n,\phi,k}^g - \mathbb{E}[\nu_{n,\phi,k}^g]$ be the centered version of the measure $\nu_{n,\phi,k}^g$.

Henceforth, we call $N_{n,\phi,k}^g$ the "$k$-nearest neighbors spacing statistic," or "empirical nearest neighbor $\phi$-divergence"; it provides a high-dimensional analog of the statistic (1.1). Our main concern is with the limit theory of $\nu_{n,\phi,k}^g$ and $N_{n,\phi,k}^g$.

The statistic $N_{n,\phi,k}^g$ provides an empirical measure of the discrepancy of the proposed distribution $G$ from the (typically unknown) true distribution $F$. For example, if $k=1$, then equating $D^f_{i,n,1}$ with its approximate expected value of $1/n$ yields the approximation $N^g_{n,\phi,1} \approx \sum_i \phi(D^g_{i,n,1}/D^f_{i,n,1})$, and thus $N^g_{n,\phi,1}$ provides a naive empirical estimate for the so-called $\phi$-*divergence* [1, 2, 3, 9, 10, 11] of $g$ from $f$ which is defined by

$$(2.4) \qquad I_\phi(g,f) := \int_A \phi\left(\frac{g(x)}{f(x)}\right) f(x)\,dx.$$

In general, $I_\phi(g,f)$ is possibly negative, and $I_\phi(g,f) = I_{\phi^*}(f,g)$ where $\phi^*(x) := x\phi(x^{-1})$. Also,

$$(2.5) \qquad I_\phi(f,f) = \phi(1);$$

$$(2.6) \qquad I_\phi(g,f) \geq I_\phi(f,f) \quad \text{if } \phi \text{ is convex.}$$

Choices of $\phi \in \mathcal{F}$ figuring prominently in estimation and decision theory include:

- $\phi_0(x) := -\log x$ defines Kullback–Leibler information (also called the modified log-likelihood ratio statistic or relative entropy) and is used in maximum spacing methods,
- $\phi_{1/2}(x) := 2(1-\sqrt{x})^2$ yields the square of the Hellinger distance,
- $\phi_1(x) := x\log x$ yields the log-likelihood ratio statistic or I-divergence of Kullback–Leibler,
- $\phi_2(x) := (x-1)^2/2$ yields the chi-squared divergence, and
- $\phi^{(r)}(x) := x^r$ yields information gain of order $r, r > 0$.

The $\phi$-divergences $N^g_{n,\phi,k}$ and $I_\phi(g,f)$ ("coefficients of divergence") are used heavily in goodness-of-fit tests [36] and are useful in characterizing the amount of information of one distribution contained in another [36, 37]. Nearest neighbor cells have been used in goodness-of-fit tests in multidimensions in [6, 38, 43], among others. Note that (2.6) shows $I_{\phi_0}(g,f)$ and $I_{\phi_1}(g,f)$ are nonnegative, and that $I_{\phi^{(1/2)}}(g,f)$ is symmetric in $f$ and $g$.

The following integral transforms of $\phi$ (defined for $\beta > 0$) arise naturally in the asymptotic analysis of $\nu^g_{n,\phi,k}$ (the random variables $\Gamma_k$ were defined in Section 2.1):

$$(2.7) \quad M_{\phi,k}(\beta) := \mathbb{E}[\phi(\beta\Gamma_k)],$$

$$(2.8) \quad \Delta_{\phi,k}(\beta) := (k+1)M_{\phi,k}(\beta) - kM_{\phi,k+1}(\beta),$$

$$V_{\phi,k}(\beta) := M_{\phi^2,k}(\beta)$$

$$(2.9) \qquad\qquad + \int_{\mathbb{R}^d} [\mathbb{E}[\phi(\beta|C_k(\mathbf{0},\mathcal{H}\cup y)|)\phi(\beta|C_k(y,\mathcal{H}\cup\mathbf{0})|)]$$

$$- M_{\phi,k}(\beta)^2]\,dy.$$

Note that $M_{\phi,1}(x) = (1/x)\hat{\phi}(1/x)$, where $\hat{\phi}$ denotes the Laplace transform of $\phi$.

2.3. *A general CLT for $\phi$-divergences.* The following general central limit theorem, our main result, establishes convergence of $n^{-1/2}\langle h, \overline{\nu}^g_{n,\phi,k}\rangle$ to a mean zero normal random variable whose variance is a weighted average of the functions $V_{\phi,k}$ and $\Delta_{\phi,k}$. For $h \in \mathcal{B}(A)$, we define the *h-weighted $\phi$-divergence* of $g$ from $f$ by

$$I_\phi(g,f,h) := \int_A f(x)\phi\left(\frac{g(x)}{f(x)}\right)h(x)\,dx,$$

which in the case $h \equiv 1$ reduces to the $\phi$-divergence $I_\phi(f,g)$ defined at (2.4). Also, for $h, h_1, h_2$ in $\mathcal{B}(A)$ and $\phi \in \mathcal{F}$, we define the functions $h^2, h_1h_2, \phi^2$ pointwise, that is, $h^2(x) = (h(x))^2$ and so on.

In the theorem below, since the formula (2.10) is rather concise, we expand it in (2.11). We prove the theorem in Section 4, referring to [4] for some of the details.

THEOREM 2.1. *Suppose that either $\phi \in \mathcal{F}_0, d=1$, or $\mathcal{K} = \mathbb{R}^d$. As $n \to \infty$, it is the case that for $h \in \mathcal{B}(A)$,*

$$(2.10) \qquad n^{-1}\operatorname{Var}[\langle h, \nu^g_{n,\phi,k}\rangle] \to I_{V_{\phi,k}}(g,f,h^2) - (I_{\Delta_{\phi,k}}(g,f,h))^2$$

$$(2.11) \qquad = \int_A h^2(x) V_{\phi,k}\left(\frac{g(x)}{f(x)}\right) f(x)\,dx$$
$$- \left(\int_A h(x) \Delta_{\phi,k}\left(\frac{g(x)}{f(x)}\right) f(x)\,dx\right)^2$$

*and*

$$(2.12) \qquad n^{-1/2}\langle h, \overline{\nu}^g_{n,\phi,k}\rangle \xrightarrow{\mathcal{D}} N(0, I_{V_{\phi,k}}(g,f,h^2) - (I_{\Delta_{\phi,k}}(g,f,h))^2).$$

Putting $h \equiv 1$ in Theorem 2.1 yields a CLT for the empirical $\phi$-divergence $N^g_{n,\phi,k}$:

$$n^{-1/2}(N^g_{n,\phi,k} - \mathbb{E}N^g_{n,\phi,k}) \xrightarrow{\mathcal{D}} N(0, I_{V_{\phi,k}}(g,f) - (I_{\Delta_{\phi,k}}(g,f))^2).$$

For practical purposes, it is of use to compute $V_{\phi,k}$, and the next two results show how to simplify the expression (2.9) in some special cases. Using these simplifications, we may explicitly identify $V_{\phi,k}$ for certain choices of $\phi$, as shown in Section 3.

The first of our simplifications applies when $\mathcal{K} \neq \mathbb{R}^d$, and either $k = 1$ or $d = 1$. The latter case is particularly relevant to the study of spacings (see Section 2.4).

PROPOSITION 2.1. *If $\mathcal{K} \neq \mathbb{R}^d$, and either $d = 1$ or $k = 1$, then for all $\beta > 0$, we have*

$$
\begin{aligned}
V_{\phi,k}(\beta) = {} & M_{\phi^2,k}(\beta) + 2kM_{\phi,k}(\beta)(M_{\phi,k}(\beta) - M_{\phi,k+1}(\beta)) \\
& + 2\sum_{j=1}^{k-1} \mathrm{Cov}[\phi(\beta\Gamma_k), \phi(\beta(\Gamma_{k+j} - \Gamma_j))],
\end{aligned}
\tag{2.13}
$$

*the sum being interpreted as zero for $k = 1$.*

Our second simplifying formula for $V_{\phi,k}$ is applicable when $k = 1$, $\mathcal{K} = \mathbb{R}^d$, and $\phi$ is differentiable with $\lim_{t\downarrow 0} \phi(t) = 0$. This will provide limiting distributions for some cases of interest, including information gain and log-likelihood in high dimensions (see Section 3.2). For $s, t, u \in \mathbb{R}^+$, let $I(s,t,u)$ be the volume of the intersection of two balls in $\mathbb{R}^d$, with respective volumes $s$ and $t$, at a distance $u$ apart. Set

$$
J_d(s,t) := \int_{\max(s,t)}^{\infty} [e^{I(s,t,(u/\omega_d)^{1/d})} - 1]\, du.
\tag{2.14}
$$

PROPOSITION 2.2. *Suppose that $\mathcal{K} = \mathbb{R}^d$ and that $\phi \in \mathcal{F}$ is differentiable with $\lim_{t\downarrow 0} \phi(t) = 0$. Then for all $\beta > 0$,*

$$
V_{\phi,1}(\beta) = M_{\phi^2,1}(\beta) + \beta^2 \int_0^\infty \int_0^\infty \phi'(\beta s)\phi'(\beta t) e^{-(s+t)}[J_d(s,t) - \max(s,t)]\, ds\, dt
$$

*provided that the integral exists.*

REMARKS. (i) (*Related work*) Bickel and Breiman [6], and subsequently Schilling [38], consider the functionals $N_{n,\phi,1}^g$ when $\phi(x) = \exp(-x)$ and $\mathcal{K} = \mathbb{R}^d$. Using the approximation $D_{i,n,1}^g \approx g(X_i)|C_1(X_i, \mathcal{X}_n)|$, they establish a CLT for the empirical process of nearest neighbor distances, but do not consider convergence of the associated random measures. Zhou and Jamalamadaka [43] establish the central limit theory for the functionals $N_{n,\phi,1}^g$ for certain $\phi$ of bounded variation for the case $g = f$ as well as for the case involving a sequence of appropriately converging alternatives. Strong limit theorems for multivariate spacings using general "shapes" are given by Deheuvels et al. [14].

(ii) (*Finite-dimensional CLT*) By standard arguments based on the Cramér–Wold device, it is straightforward to deduce from Theorem 2.1 the convergence of the finite-dimensional distributions of $n^{-1/2}\overline{\nu}_{n,\phi,k}^g$ as $n \to \infty$ [i.e., the convergence of the $m$-vector $n^{-1/2}(\langle h_1, \overline{\nu}_{n,\phi,k}^g\rangle, \ldots, \langle h_m, \overline{\nu}_{n,\phi,k}^g\rangle)$ for all $h_1, \ldots, h_m$ in $\mathcal{B}(A)$] to those of a mean zero finitely additive Gaussian field with covariance kernel

$$
(h_1, h_2) \mapsto I_{V_{\phi,k}}(g, f, h_1 h_2) - I_{\Delta_{\phi,k}}(g, f, h_1) I_{\Delta_{\phi,k}}(g, f, h_2).
\tag{2.15}
$$

(iii) (*Poisson CLT*) For $\lambda > 0, k \in \mathbb{N}$, the Poisson analog of measure $\nu_{n,\phi,k}^g$ is

(2.16) $$\mu_{\lambda,\phi,k}^g := \sum_{x \in \mathcal{P}_\lambda} \phi(\lambda G(C_k(x, \mathcal{P}_\lambda))) \delta_x,$$

and its total measure is a Poissonized version of $N_{n,\phi,k}^g$. Our approach yields a proof (see Proposition 4.1 below) that if $\phi \in \mathcal{F}$ and $h \in \mathcal{B}(A)$, then as $\lambda \to \infty$,

(2.17) $$\lambda^{-1} \operatorname{Var}[\langle h, \mu_{\lambda,\phi,k}^g \rangle] \to I_{V_{\phi,k}}(g, f, h^2)$$

and $\lambda^{-1/2} \overline{\mu}_{\lambda,\phi,k}^g$ converges in law to a mean zero Gaussian field with covariance kernel $(h_1, h_2) \mapsto I_{V_{\phi,k}}(g, f, h_1 h_2)$ (here $\overline{\mu}_{\lambda,\phi,k}^g := \mu_{\lambda,\phi,k}^g - \mathbb{E}[\mu_{\lambda,\phi,k}^g]$).

(iv) (*Law of large numbers, limits are distribution free*) Our approach (see also [31]) also yields a weak law of large numbers, namely

$$n^{-1} \langle h, \nu_{n,\phi,k}^g \rangle \xrightarrow{L^2} I_{M_{\phi,k}}(g, f, h) \qquad \forall h \in \mathcal{B}(A), \phi \in \mathcal{F}.$$

By taking $h \equiv 1$, we obtain a weak law of large numbers for the $k$-nearest neighbors spacing statistic $N_{n,\phi,k}^g$. Combining this with Theorem 2.1 and taking $g = f$, we see from (2.5) that the limiting mean of $n^{-1} \langle h, \nu_{n,\phi,k}^f \rangle$ and the limiting variance and distribution of $n^{-1/2} \langle h, \overline{\nu}_{n,\phi,k}^f \rangle$ do not depend on $f$ for $h \equiv 1$ (and, in fact, for any $h$). Therefore, the nearest neighbor functionals are *asymptotically distribution free under the null hypothesis* $g = f$ and have asymptotic variance $V_{\phi,k}(1) - (\Delta_{\phi,k}(1))^2$. A possible goodness-of-fit test would be to take the density $g$ to be tested, compute the functional $N_{n,\phi,1}^g$ and see whether the cumulative distribution function is close to the $N(0, V_{\phi,1}(1) - (\Delta_{\phi,1}(1))^2)$ cumulative distribution function.

(v) (*Voronoi cells*) Volumes of nearest neighbor cells are computationally attractive and have correlations decaying exponentially with the distance between cell centers. Defining point measures analogous to (2.3) based on cells generated by any locally defined Euclidean graph (e.g., Voronoi cells) leads to similar CLTs, adding to the laws of large numbers given in [24].

(vi) (*Properties of limiting variance*) In most of our examples, $\Delta_{\phi,k}$ is strictly positive, showing that Poissonization leads to a larger limiting variance. When $V_{\phi,k}$ is convex, which is the case when $k = 1$, $\phi(x) = x^r, r \in [1, \infty)$ or when $\phi(x) = x \log x$ (see Section 3.1), then inequality (2.6) implies that the limiting variance over Poisson samples is minimized when $g = f$.

2.4. *Asymptotic normality of sum functions of spacings*. In dimension $d = 1$, if $g$ is a probability density with distribution function $G$ on $[c_1, c_2]$, then the generalization to $k$-spacings of the empirical $\phi$-divergence defined at (1.1) is the classical $k$-spacing statistic defined by

$$(2.18) \qquad S_{n,\phi,k}^g := \sum_{i=1}^{n-k} \phi(n[G(X_{(i+k)}) - G(X_{(i)})]).$$

Developing the limit theory for $S_{n,\phi,k}^g$ over continuous samples is important in goodness-of-fit tests. We can apply the general theory of Section 2.3 by putting $d = 1$ and $\mathcal{K} = (0, \infty)$. Then the width of $C_k^{\mathcal{K}}(x, \mathcal{X})$ is the distance between $x$ and its $k$th nearest neighbor in $\mathcal{X}$ "to the right." Thus, the $k$-nearest neighbors spacing statistic $N_{n,\phi,k}^g$, defined by (2.3), is the same as $S_{n,\phi,k}^g$ but with the sum in (2.18) extended to $n$ terms and with $X_{(j)} := c_2$ if $j > n$.

To better match the existing literature, we consider a modified version of Theorem 2.1 in which we redefine $C_k^{\mathcal{K}}(x, \mathcal{X})$ to be the empty set whenever $\text{card}(\mathcal{X} \cap (x + \mathcal{K}) \setminus x) < k$, and set $\phi(0) = 0$. Denote by $\nu_{n,\phi,k}^*$ the analog of $\nu_{n,\phi,k}^g$ under this modification (here we suppress the dependence on $g$), that is,

$$(2.19) \qquad \nu_{n,\phi,k}^* := \sum_{i=1}^{n-k} \phi(n[G(X_{(i+k)}) - G(X_{(i)})]) \delta_{X_i}.$$

The corresponding centered measure is then denoted $\bar{\nu}_{n,\phi,k}^*$. If $d = 1$ and $\mathcal{K} = (0, \infty)$, the total measure of $\nu_{n,\phi,k}^*$ is indeed equal to $S_{n,\phi,k}^g$.

THEOREM 2.2 (Gaussian limit for sum functions of spacings). *Let $A := [c_1, c_2]$, $\mathcal{K} = (0, \infty)$ and $\phi \in \mathcal{F}$. Then the conclusion of Theorem 2.1 holds with $\nu_{n,\phi,k}^g$ replaced by $\nu_{n,\phi,k}^*$. Moreover, in this case, $V_{\phi,k}(\beta)$ is given by (2.13).*

The proof of Theorem 2.2 is a straightforward modification of the proof of Theorem 2.1; see Theorem 2.2 of [4] for details.

Applications of Theorem 2.2 are given in Section 3. This result, like our main result, shows that sum functions of spacings are asymptotically distribution free under the null hypothesis $f = g$.

REMARKS. (i) Darling [13] undertook the first systematic study of the functionals $S_{n,\phi,k}$ when $k = 1$, but restricted attention to uniform samples. Theorem 2.2 generalizes Holst [21], as well as earlier work of Cressie [8], who proves asymptotic normality (but not convergence of $\nu_{n,\phi,k}^*$ against bounded test functions) for sum functions of $k$-spacings over *uniform* points. Holst uses a generalization of LeCam's method and a CLT for $k$-dependent random variables. In $d = 1$, Holst and Rao [22] prove asymptotic normality of $S_{n,\phi,k}^g$ under "somewhat stringent conditions" on $f$ and $g$. Mirakhmedov [28] considers the error term in the CLT for the functionals $S_{n,\phi,k}$ when $F$ and $G$ are the uniform distribution on $[0,1]$. For nonuniform samples, the asymptotics of $S_{n,\phi,k}^g$ have been widely studied under the assumption that $G$ runs through a sequence of alternatives $G_n$ approaching the uniform distribution; see Hall [20], Kuo and Rao [27] and del Pino [34]. Khashimov [26] establishes asymptotic normality of $S_{n,\phi,k}^1$ under rather technical differentiability conditions on $\phi$ and $f$.

(ii) The approach used here also yields a weak law of large numbers, namely convergence in mean-square of $n^{-1}S_{n,\phi,k}^g$ to $I_{M_\phi}(g,f)$. This extends the corresponding weak laws in [25]; see also [17]. Analogous results hold for nonoverlapping $k$-spacings [39].

2.5. *Divergences based on cells of fixed radius.* Instead of considering point measures based on spacings, we now consider using cells of *fixed* radius depending on a continuous $g: A \to \mathbb{R}^+$ and a parameter $t$. Thus, given $\phi \in \mathcal{F}$ and $t > 0$, we define

$$H_{n,\phi}^{g,t} := \tfrac{1}{2} \sum_{x \in \mathcal{X}_n} \phi(\operatorname{card}\{\mathcal{X}_n \cap B_{t(ng(x))^{-1/d}}(x)\} - 1).$$

When $\phi(x) \equiv x$ and $g \equiv 1$, then $H_{n,\phi}^{g,t}$ counts the total number of pairs of points in $\mathcal{X}_n$ distant at most $n^{-1/d}t$ from each other.

The following CLT is obtained by modifying the proof of Theorem 2.1; we refer to Theorem 2.3 of [4] for details.

THEOREM 2.3 ([4]). *(Gaussian limit for the number of pairs of points within distance $t$).* For all continuous $g: A \to \mathbb{R}^+, t > 0$, and $\phi \in \mathcal{F}$, there is a constant $\sigma_{t,\phi,g}^2(f)$ such that as $n \to \infty$ we have $n^{-1}\operatorname{Var}[H_{n,\phi}^{g,t}] \to \sigma_{t,\phi,g}^2(f)$ and

$$n^{-1/2}(H_{n,\phi}^{g,t} - \mathbb{E}H_{n,\phi}^{g,t}) \xrightarrow{\mathcal{D}} N(0, \sigma_{t,\phi,g}^2(f)).$$

REMARKS. The limiting variance $\sigma_{t,\phi,g}^2(f)$ takes the form of the right-hand side of (2.10) with $h \equiv 1$ and with the functions $V_{\phi,k}$ and $\delta_{\phi,k}$ suitably modified; see [4] for details.

Various authors have studied $H_{n,\phi}^{g,t}$ when $\phi(x) \equiv x$ and $g \equiv 1$; see Chapter 3 of [30] and references therein. Jammalamadaka and Zhou [23] and also L'Écuyer et al. [16] consider $H_{n,\phi}^{g,t}$ from the point of view of goodness-of-fit tests, (the latter reference only for uniform samples). Penrose [30] (Chapter 4) proves that the finite-dimensional distributions of the process $H_{n,\phi}^{g,t}, t > 0$, converge to those of a Gaussian process.

**3. Applications.**

3.1. *Classical spacing statistics.* For many tests involving goodness-of-fit (Dudewicz et al. [15], Blumenthal [7], Cressie [8], Holst and Rao [22], del Pino [34], Weiss [42]) and parametric estimation (Ghosh and Jammalamadaka [19]), it is important to know the asymptotic distribution of $S^g_{n,\phi,k}$ [defined at (2.18)] for arbitrary $g$ and $f$ and for various choices of $\phi$. The following provides some illustrative examples. For simplicity of exposition, we have chosen to state our central limit theorems for the statistic $S^g_{n,\phi,k}$; the results for associated random point measures (2.19) are given in [4]. Throughout Section 3.1, we write $V^S_{\phi,k}$ for the value of $V_{\phi,k}$ given by (2.13).

3.1.1. *Limit theory for logarithms of spacings.* Let

$$S^g_{n,\log,k} := \sum_{i=1}^{n-k} \log(n[G(X_{(i+k)}) - G(X_{(i)})])$$

denote the sum of the logarithmic $k$-spacings. Setting $\phi(x) = \log x$ in Theorem 2.2 and appealing to (2.8) and (2.13), we find a CLT for logarithms of $k$-spacings as follows.

Let $\psi$ be the di-gamma function with $\psi(k) := \sum_{i=1}^{k-1} i^{-1} - \gamma$, where $\gamma$ is Euler's constant, and let $\psi'(k) := -\sum_{i=1}^{k-1} i^{-2} + \pi^2/6$.

By Cressie [8] and Holst [21],

$$\sum_{j=1}^{k-1} \mathrm{Cov}(\log \Gamma_k, \log(\Gamma_{k+j} - \Gamma_j)) = k(k-1)\psi'(k) - (k-1).$$

Also, $\mathbb{E}[\log \Gamma_k] = \psi(k)$, so we have $2k\mathbb{E}[\log \Gamma_k](\mathbb{E}\log \Gamma_k - \mathbb{E}\log \Gamma_{k+1}) = -2\psi(k)$. Also, $\mathbb{E}[\log^2 \Gamma_k] = \psi'(k) + (\psi(k))^2$. So, combining terms and using (2.13) for $\phi(x) = \log x$ gives

$$(3.1) \quad V^S_{\log,k}(1) = \psi'(k) + (\psi(k))^2 - 2\psi(k) + 2[k(k-1)\psi'(k) - (k-1)].$$

By (2.8), we have $\Delta_{\log,k}(1) = (k+1)\psi(k) - k\psi(k+1) = \psi(k) - 1$.

Using simple relations such as $\mathrm{Cov}(\log \beta X, \log \beta Y) = \mathrm{Cov}(\log X, \log Y)$, it is straightforward to deduce that $V^S_{\log,k}(\beta) = V^S_{\log,k}(1) + \log^2 \beta + 2\log \beta(\psi(k) - 1)$ and $\Delta_{\log,k}(\beta) = \Delta_{\log,k}(1) + \log \beta$. Substituting this into Theorem 2.2, putting $\tau_k := (2k^2 - 2k + 1)\psi'(k) - 2k + 1$, and rearranging terms yields the following corollary.

COROLLARY 3.1 (CLT for logarithmic $k$-spacings). *Let $X, X_1, X_2, \ldots$ be i.i.d. with density $f$ on $[0,1]$. As $n \to \infty$, $n^{-1} \mathrm{Var}[S^g_{n,\log,k}] \to \tau_k + \mathrm{Var}[\log(\frac{f(X)}{g(X)})]$ and*

$$n^{-1/2}(S^g_{n,\log,k} - \mathbb{E}S^g_{n,\log,k}) \xrightarrow{\mathcal{D}} N\left(0, \tau_k + \mathrm{Var}\left[\log\left(\frac{f(X)}{g(X)}\right)\right]\right).$$

REMARKS. When $A = [0, 1]$ and $f \equiv g \equiv 1$, the CLT for $S_{n,\log,k}^g$ was established by Darling (Section 7 of [13]) for $k = 1$ and later by Holst [21] and Cressie [8] for general $k$. When the $X_i$ have a step density, Cressie shows asymptotic normality of $S_{n,\log,k}^g$ including cases when $k \to \infty$. Czekała (Theorem 1 of [12]) apparently rediscovered Cressie's result. Shao and Hahn [40] treat general densities for $k = 1$, although their proof depends upon interchanging limits in order to pass from step densities to arbitrary densities. When $k = 1$, Blumenthal (Theorem 2 of [7]), proves Corollary 3.1 for densities $f$ satisfying special conditions. Corollary 3.1 extends all of these results to $f$ and $g$ bounded away from zero and infinity, resolving a conjecture of Darling ([13], page 249) affirmatively.

3.1.2. *Information gain of order $r$.* Let $\phi(x) = x^r, r > 0$. We write $S_{n,r,1}^g$ to denote $S_{n,\phi,1}^g$, also known as Rényi's information gain (I-divergence) of order $r$ in $d = 1$, that is,

$$S_{n,r,1}^g := \sum_{i=1}^{n-1} (n[G(X_{(i+1)}) - G(X_{(i)})])^r.$$

Let $w_r := -2r\Gamma^2(r+1) + \Gamma(2r+1)$ and $t_r := \Gamma(r+1)(1-r)$. It is a simple matter to verify via (2.13) and (2.8), respectively, that for all $\beta > 0$,

$$V_{\phi,1}^S(\beta) := w_r \beta^{2r} \quad \text{and} \quad \Delta_{\phi,1}(\beta) := 2\mathbb{E}[\phi(\beta\Gamma_1)] - \mathbb{E}[\phi(\beta\Gamma_2)] = t_r \beta^r.$$

Put

$$\sigma_r^2(f, g) := w_r \int_A \left(\frac{g(x)}{f(x)}\right)^{2r} f(x) \, dx - t_r^2 \left(\int_A \left(\frac{g(x)}{f(x)}\right)^r f(x) \, dx\right)^2.$$

Theorem 2.2 yields the following corollary.

COROLLARY 3.2 (Gaussian limits for information gain). *Let $X_1, X_2, \ldots$ be i.i.d. with density $f$ on $A := [c_1, c_2]$. As $n \to \infty$, we have for all $h \in \mathcal{B}(A)$*

$$n^{-1} \operatorname{Var}[S_{n,r,1}^g] \to \sigma_r^2(f, g)$$

*and* $n^{-1/2}(S_{n,r,1}^g - \mathbb{E} S_{n,r,1}^g) \xrightarrow{\mathcal{D}} N(0, \sigma_r^2(f, g))$.

REMARKS. It is easy to verify using [5] that $\sigma_r^2(f, g) > 0$ except when $r = 1$. Corollary 3.2 extends upon the CLTs of Darling [13] (uniform case) and Weiss [42]. Moran [29] proved a CLT for the functional $S_{n,1,1}^g$ over uniform random variables.

3.1.3. *Limit theory for log-likelihood ratio.* Let $\phi(x) = x \log x$. Consider the log-likelihood point measure

$$\nu_{n,\phi,1}^g := \sum_{i=1}^{n-1} \phi(n[G(X_{(i+1)}) - G(X_{(i)})])\delta_{X_i}$$

and let $S_\phi^g$ denote the total mass of this measure, also called the log-likelihood statistic. Again, denoting Euler's constant by $\gamma$, we have for $\beta > 0$ that

$$\mathbb{E}[\beta\Gamma_1 \log(\beta\Gamma_1)] = \beta \log \beta + \beta(1 - \gamma);$$
$$\mathbb{E}[\beta\Gamma_2 \log(\beta\Gamma_2)] = 2\beta \log \beta + \beta(3 - 2\gamma);$$
$$\mathbb{E}[(\beta\Gamma_1 \log(\beta\Gamma_1))^2] = \beta^2[2(\log \beta)^2 + (6 - 4\gamma)\log \beta + 2 + \pi^2/3 - 6\gamma + 2\gamma^2].$$

Using these in (2.13) and (2.8), respectively, it is easily verified that

$$V_{\phi,1}^S(\beta) := \left(\frac{\pi^2}{3} - 2\right)\beta^2 \quad \text{and} \quad \Delta_{\phi,1}(\beta) := 2\mathbb{E}\phi(\beta\Gamma_1) - \mathbb{E}\phi(\beta\Gamma_2) = -\beta.$$

Put

$$\sigma_\phi^2(f,g) := \left(\frac{\pi^2}{3} - 2\right) \int_A \frac{g^2(x)}{f(x)}\,dx - \left(\int_A g(x)\,dx\right)^2.$$

Let $X$ have density $f$ and note that since $g$ is a density we have

$$\sigma_\phi^2(f,g) = \left(\frac{\pi^2}{3} - 2\right) \text{Var}\left[\frac{g(X)}{f(X)}\right] + \frac{\pi^2}{3} - 3.$$

Using the above values for $V_{\phi,1}$, $\Delta_{\phi,1}$, $\sigma_\phi^2(f,g)$, and applying Theorem 2.2 for $\phi(x) = x \log x$ yields the following corollary.

COROLLARY 3.3 (Gaussian limit for log-likelihood). *Let $X_1, X_2, \ldots$ be i.i.d. with density $f$ on $A := [c_1, c_2]$. As $n \to \infty$, $n^{-1}\text{Var}[S_{n,\phi,1}^g] \to \sigma_\phi^2(f,g)$ and*

$$n^{-1/2}(S_{n,\phi,1}^g - \mathbb{E}S_{n,\phi,1}^g) \xrightarrow{\mathcal{D}} N(0, \sigma_\phi^2(f,g)).$$

REMARKS. Corollary 3.3 extends the results of Gebert and Kale [18], who assume uniformity of $X_i$ and Czekała (Theorem 2 of [12]), who assumes that $X_i$ have a step density. van Es [41] establishes asymptotic normality for $S_\phi^g$ whenever $k, n \to \infty$, $k = o(n^{1/2})$, and $f: A \to [0, \infty)$ is Lipschitz.

3.2. *Information gain and log-likelihood in high dimensions.* In this section, we put $k = 1$ and $\mathcal{K} = \mathbb{R}^d$.

3.2.1. *Information gain of order $r$.* Let $\phi(x) = x^r, r \in \mathbb{R}^+$, so that $N_{n,\phi,1}^g$ defined by (2.3) yields Rényi's information gain (I-divergence) of order $r$. For all $r \in \mathbb{R}^+$, define the constant

$$K_r := r^2 \int_0^\infty \int_0^\infty s^{r-1} t^{r-1} e^{-(s+t)} [J_d(s,t) - \max(s,t)] \, ds \, dt,$$

with $J_d(s,t)$ given by (2.14). Since $\phi$ satisfies the conditions of Proposition 2.2 and since $\mathbb{E}[\phi^2(\Gamma_1)] = \Gamma(2r+1)$, the following is immediate.

LEMMA 3.1. *For all $\beta > 0$ and for $\phi(x) = x^r, r > 0$, we have that $V_{\phi,1}(\beta) = \beta^{2r}[\Gamma(2r+1) + K_r]$.*

Note that $\beta^{2r} = \phi^2(\beta)$. Combining Lemma 3.1 with Theorem 2.1 yields the following CLT for $N_{n,\phi,1}^g$.

COROLLARY 3.4. *Let $\phi(x) = x^r, r > 0$. Then as $n \to \infty$ $n^{-1} \operatorname{Var}[N_{n,\phi,1}^g]$ converges to $[\Gamma(2r+1) + K_r]I_{\phi^2,1}(g,f) - (I_{\Delta_{\phi,1}}(f,g))^2$, and*

$$n^{-1/2}(N_{n,\phi,1}^g - \mathbb{E} N_{n,\phi,1}^g)$$
$$\xrightarrow{\mathcal{D}} N(0, [\Gamma(2r+1) + K_r]I_{\phi^2,1}(g,f) - (I_{\Delta_{\phi,1}}(g,f))^2).$$

3.2.2. *Log-likelihood.* When $\phi(x) = x \log x$, $N_{n,\phi}^g$ defined by (2.3) yields the log-likelihood statistic. To apply Theorem 2.1, we define the constants

$$I_1 := \int_0^\infty \int_0^\infty (\log s + 1)(\log t + 1) e^{-(s+t)} [J_d(s,t) - \max(s,t)] \, ds \, dt,$$
$$I_2 := \int_0^\infty \int_0^\infty (\log s + 1) e^{-(s+t)} [J_d(s,t) - \max(s,t)] \, ds \, dt$$

and

$$I_3 := \int_0^\infty \int_0^\infty e^{-(s+t)} [J_d(s,t) - \max(s,t)] \, ds \, dt.$$

Also, set $K_1 := 2 + \frac{\pi}{3} - 6\gamma + 2\gamma^2 + I_1$, $K_2 := 6 - 4\gamma + 2I_2$ and $K_3 := 2 + I_3$. The following is an easy consequence of Proposition 2.2.

LEMMA 3.2. *For $\phi(x) = x \log x$, $V_{\phi,1}(\beta) = \beta^2(K_1 + K_2 \log \beta + K_3(\log \beta)^2)$, $\beta > 0$.*

Theorem 2.1 yields a CLT for the log-likelihood functional $N_{n,\phi,1}^g$. Put

$$\sigma_\phi^2(f,g) := \int_A \left(\frac{g(x)}{f(x)}\right)^2 \left[K_1 + K_2 \log\left(\frac{g(x)}{f(x)}\right) + K_3 \left(\log \frac{g(x)}{f(x)}\right)^2\right] f(x) \, dx.$$

COROLLARY 3.5. *Let $\phi(x) = x \log x$. Then as $n \to \infty$, $n^{-1} \operatorname{Var}[N^g_{n,\phi,1}] \to \sigma^2_\phi(f,g) - (I_{\Delta_{\phi,1}}(g,f))^2$ and*

$$n^{-1/2}(N^g_{n,\phi,1} - \mathbb{E} N^g_{n,\phi,1}) \xrightarrow{\mathcal{D}} N(0, \sigma^2_\phi(f,g) - (I_{\Delta_{\phi,1}}(g,f))^2).$$

**4. Proof of Theorem 2.1.** The proof of Theorem 2.1 involves expressing the Poissonized measure $\mu^g_{\lambda,\phi,k}$ [see (2.16)] as a sum of weakly spatially dependent terms, allowing us to establish convergence of the variance of the measure $\mu^g_{\lambda,\phi,k}$ (Proposition 4.1). Although the measures in question share neither the same representation nor the same scaling properties as those considered in previous work [5, 31, 32, 33], once we have shown the crucial variance convergence for measures defined in terms of Poisson samples, we can draw upon some well-established dependency graph techniques [31, 32, 33] to deduce a Poissonized version of Theorem 2.1. Using arguments in [32], we may de-Poissonize and deduce Theorem 2.1 when $\phi$ is bounded on $(0, 1]$. Deducing Theorem 2.1 for general $\phi$ requires extra technical effort.

Recall that for all $a > 0$, $\mathcal{H}_a$ is a homogeneous Poisson point process of intensity $a$ on $\mathbb{R}^d$. Suppose we fix the set $A \subset \mathbb{R}^d$, the densities $f$ and $g$ and their corresponding distributions $F$ and $G$ on $\mathbb{R}^d$, as described in Section 2.1.

For all $\lambda > 0, x \in \mathbb{R}^d$, and all finite $\mathcal{X} \subset \mathbb{R}^d$, we lighten the notation and write $C(x, \mathcal{X})$ for $C^{\mathcal{K}}_k(x, \mathcal{X})$ given by (2.2). For all Borel $B \subset \mathbb{R}^d$, define the numbers $\Phi_\lambda(x, \mathcal{X}) = \Phi^g_\lambda(x, \mathcal{X})$ and $\xi_\lambda(x, \mathcal{X}, B) := \xi^g_\lambda(x, \mathcal{X}, B)$ by

$$\Phi_\lambda(x, \mathcal{X}) := \phi(\lambda G(C(x, \mathcal{X})); \qquad \xi_\lambda(x, \mathcal{X}, B) := \Phi_\lambda(x, \mathcal{X}) \delta_x(B).$$

Recalling that $\mathcal{X}_n := \{X_1, \ldots, X_n\}$, we have

$$\nu^g_{n,\phi,k}(B) = \sum_{i=1}^n \xi^g_n(X_i, \mathcal{X}_n, B).$$

The (signed) point measure $\xi_\lambda(x, \mathcal{X}, \cdot)$ is determined by $x$ and $\mathcal{X}$, in a similar manner to the measures considered in [32], but here, unlike in [32], the measure $\xi_\lambda(x, \mathcal{X}, \cdot)$ is not obtained by scaling the measure $\xi_1$ (because the function $g$ enters in a more complicated way into the definition of $\xi_\lambda$ here) so we cannot directly apply results from [32]. We write $\langle h, \xi_\lambda(x, \mathcal{X}) \rangle$ for $\int_{\mathbb{R}^d} h(y) \xi_\lambda(x, \mathcal{X}, dy)$.

For locally finite $\mathcal{X} \subset \mathbb{R}^d$ and $x \in A$, we define

(4.1) $\qquad \xi^{g,x}_\infty(\mathcal{X}) := \xi^{g,x,k}_\infty(\mathcal{X}) := \phi(g(x)|C(\mathbf{0}, \mathcal{X})|).$

For all $x \in A$ and given $k \in \mathbb{N}$, let $t_0(x)$ denote the infimum of all $t$ with the property that $B^{\mathcal{K}}_u(x) \cap A$ is the same for all $u \geq t$. Given also a locally finite set $\mathcal{X} \subset A$ and $\mathcal{K}$, and writing $\#(\cdot)$ for card$(\cdot)$ and $\mathcal{X} \setminus x$ for $\mathcal{X} \setminus \{x\}$ here, define

$$R_1(x, \mathcal{X})$$
$$:= \begin{cases} \inf\{t \in \mathbb{R}^+ : \#(B_t^{\mathcal{K}}(x) \cap \mathcal{X} \setminus x) \geq k\}, & \text{if } \#((x+\mathcal{K}) \cap \mathcal{X} \setminus x) \geq k, \\ t_0(x), & \text{otherwise.} \end{cases}$$

Thus, $R_1(x, \mathcal{X})$ is the distance between $x$ and its $k$th nearest neighbor in $\mathcal{X}$ in the direction of the cone $\mathcal{K}$ or if no such neighbor exists, the furthest one has to look from $x$ to ascertain that this is the case. For $\lambda > 0$, let $R_\lambda(x, \mathcal{X}) := \lambda^{1/d} R_1(x, \mathcal{X})$. The following lemma establishes the equivalent of the "exponential stabilization" conditions discussed in [32]. For a proof, see [4].

LEMMA 4.1 ([4], Lemma 4.3). *It is the case that*

$$\limsup_{t \to \infty} \sup_{x \in A, \lambda \geq 1} t^{-1} \log P[R_\lambda(x, \mathcal{P}_\lambda) > t] < 0 \tag{4.2}$$

*and*

$$\limsup_{t \to \infty} \sup_{x \in A, \lambda \geq 1, (\lambda/2) \leq n \leq (3\lambda/2), \mathcal{A} \in \mathcal{S}_3} t^{-1} \log P[R_\lambda(x, \mathcal{X}_n \cup \mathcal{A}) > t] < 0. \tag{4.3}$$

Recall (2.16) that $\mu_{\lambda,\phi,k}^g$ denotes the Poissonized version of $\nu_{\lambda,\phi,k}^g$. The following is a Poissonized version of Theorem 2.1 for $\phi \in \mathcal{F}_0$ and is of independent interest.

PROPOSITION 4.1. *Let $\phi \in \mathcal{F}_0$ and $h \in \mathcal{B}(A)$. Then as $\lambda \to \infty$*

$$\lambda^{-1} \operatorname{Var}[\langle h, \mu_{\lambda,\phi,k}^g \rangle] \to \int_A h^2(x) V_{\phi,k}\left(\frac{g(x)}{f(x)}\right) f(x) \, dx = I_{V_{\phi,k}}(g, f, h^2) \tag{4.4}$$

*and $\lambda^{-1/2} \overline{\mu}_{\lambda,\phi,k}^g$ converges in law as $\lambda \to \infty$ to a mean zero Gaussian field with covariance kernel $(h_1, h_2) \mapsto I_{V_{\phi,k}}(g, f, h_1) I_{V_{\phi,k}}(g, f, h_2)$.*

PROOF. For simplicity, we first assume that $h$ is a.e. continuous. It is the case that

$$\begin{aligned}
\lambda^{-1} &\operatorname{Var}[\langle h, \mu_{\lambda,\phi,k}^g \rangle] \\
&= \lambda \int_A \int_A h(x) h(y) \{\mathbb{E}[\Phi_\lambda(x, \mathcal{P}_\lambda \cup y) \Phi_\lambda(y, \mathcal{P}_\lambda \cup x)] \\
&\qquad\qquad - \mathbb{E}[\Phi_\lambda(x, \mathcal{P}_\lambda)] \mathbb{E}[\Phi_\lambda(y, \mathcal{P}_\lambda)]\} f(x) f(y) \, dx \, dy \\
&\quad + \int_A h^2(x) \mathbb{E}[\Phi_\lambda^2(x, \mathcal{P}_\lambda)] f(x) \, dx.
\end{aligned} \tag{4.5}$$

We will sketch an argument (see Proposition 4.1 of [4] for details) showing that $\lambda^{-1} \operatorname{Var}[\langle h, \mu^g_{\lambda,\phi,k}\rangle]$ converges to

$$
\begin{aligned}
\int_A \int_{\mathbb{R}^d} h^2(x) & [\mathbb{E}\xi^{g,x}_\infty(\mathcal{H}_{f(x)} \cup z)\xi^{g,x}_\infty(-z + (\mathcal{H}_{f(x)} \cup \mathbf{0})) \\
& \qquad - (\mathbb{E}\xi^{g,x}_\infty(\mathcal{H}_{f(x)}))^2]f^2(x)\,dz\,dx \\
+ \int_A h^2(x) & \mathbb{E}[(\xi^{g,x}_\infty(\mathcal{H}_{f(x)}))^2]f(x)\,dx.
\end{aligned}
\tag{4.6}
$$

Putting $y = x + \lambda^{-1/d}z$ in the right-hand side in (4.5) reduces the double integral to

$$
(4.7) \quad = \int_A \int_{-\lambda^{1/d}x+\lambda^{1/d}A} h(x)h(x+\lambda^{-1/d}z)\{\cdots\}f(x)f(x+\lambda^{-1/d}z)\,dz\,dx
$$

where

$$
\begin{aligned}
\{\cdots\} := \{&\mathbb{E}[\Phi_\lambda(x, \mathcal{P}_\lambda \cup \{x+\lambda^{-1/d}z\})\Phi_\lambda(x+\lambda^{-1/d}z, \mathcal{P}_\lambda \cup x)] \\
& - \mathbb{E}[\Phi_\lambda(x, \mathcal{P}_\lambda)]\mathbb{E}[\Phi_\lambda(x+\lambda^{-1/d}z, \mathcal{P}_\lambda)]\}.
\end{aligned}
$$

By using Lemma 4.1, we may show (see [4] for details) that $\{\cdots\}$ converges to the bracketed expression in the first term of (4.6), and that the integrand in (4.7) is dominated by an integrable function of $z$ over $\mathbb{R}^d$. The convergence of the double integral in (4.5) to that in (4.6) now follows by dominated convergence, the continuity of $h$ and fourth moment bounds on $\Phi_\lambda$. To show convergence of general $h \in \mathcal{B}(A)$, we refer to [32].

Similar but easier methods show convergence of $\int_A h^2(x)\mathbb{E}[\Phi^2_\lambda(x, \mathcal{P}_\lambda)]f(x)\,dx$, completing the proof that (4.5) converges to (4.6).

For all $x \in A$, we define $V^\xi_{\phi,k}(x, 0) := 0$ and for all $a > 0$, we put

$$
\begin{aligned}
V^\xi_{\phi,k}(x, a) := & \mathbb{E}[\xi^{g,x}_\infty(\mathcal{H}_a)^2] \\
& + a\int_{\mathbb{R}^d}[\mathbb{E}\xi^{g,x}_\infty(\mathcal{H}_a \cup z)\xi^{g,x}_\infty(-z + (\mathcal{H}_a \cup \mathbf{0})) - (\mathbb{E}\xi^{g,x}_\infty(\mathcal{H}_a))^2]\,dz.
\end{aligned}
$$

Using (4.1), it is easy to see that

$$
\begin{aligned}
V^\xi_{\phi,k}(x, a) = \mathbb{E}\left[\phi\left(\frac{g(x)}{a}\Gamma_k\right)^2\right] + \int_{\mathbb{R}^d}\bigg[&\mathbb{E}\left[\phi\left(\frac{g(x)}{a}|C(\mathbf{0}, \mathcal{H} \cup y)|\right)\right. \\
& \left.\times \phi\left(\frac{g(x)}{a}|C(y, \mathcal{H} \cup \mathbf{0})|\right)\right] \\
& - \left(\mathbb{E}\left[\phi\left(\frac{g(x)}{a}\Gamma_k\right)\right]\right)^2\bigg]\,dy
\end{aligned}
\tag{4.8}
$$

and in particular, by definition of $V_{\phi,k}$ [recall (2.9)], we have

$$V_{\phi,k}^{\xi}(x, f(x)) = V_{\phi,k}\left(\frac{g(x)}{f(x)}\right).$$

By combining this with (4.6), we thus obtain the desired limiting variance (4.4).

The proof of the second part of Proposition 4.1 (i.e., convergence to the normal) follows from arguments similar to those used in Theorem 2.2 of [32], which itself follows dependency graph arguments in [33]. Here, we note that $\xi_\infty^{g,x}(\mathcal{H}_{f(x)})$ corresponds in our setting to the limiting expression from Lemma 3.4 of [32], and consequently appears in expressions for limiting variances arising from following the proofs in [32], where all expressions for limits are obtained through Lemmas 3.4 and 3.5 of [32].  □

To obtain Theorem 2.1, we shall de-Poissonize Proposition 4.1 by suitably adapting the proofs of Theorems 2.1, 2.2 and 2.3 of [32], and then extend the result to all $\phi \in \mathcal{F}$ by using truncation arguments.

PROPOSITION 4.2.  *Suppose $\phi \in \mathcal{F}_0$. Then the conclusions of Theorem 2.1 hold.*

PROOF.  Taking Proposition 4.1 as our starting point, we can follow nearly verbatim the de-Poissonization argument of Section 5 of [32] which is used there to prove Theorem 2.3 of [32]. See Proposition 4.2 of [4] for details.

To obtain the limiting variance in the present setting, it suffices to consider the corresponding limits obtained in Lemmas 3.4 and 3.5 of [32]. That is, analogously to the definition of $\delta(x,a)$ in [32], we define for all $x \in A$ and all $a > 0$

$$\Delta_{\phi,k}^{\xi}(x,a) := \mathbb{E}[\xi_\infty^{g,x}(\mathcal{H}_a)] + a \int_{\mathbb{R}^d} [\mathbb{E}\xi_\infty^{g,x}(\mathcal{H}_a \cup y) - \xi_\infty^{g,x}(\mathcal{H}_a)] \, dy$$

$$= \mathbb{E}\left[\phi\left(\frac{g(x)}{a}\Gamma_k\right)\right]$$

$$+ a \int_{\mathbb{R}^d} \mathbb{E}[\phi(g(x)|C(\mathbf{0}, \mathcal{H}_a \cup y)|) - \phi(g(x)|C(\mathbf{0}, \mathcal{H}_a)|)] \, dy.$$

Changing variables $y \to a^{1/d}y$, using (2.7) and the equivalence $a^{1/d}\mathcal{H}_a \stackrel{\mathcal{D}}{=} \mathcal{H}$, yields

$$\Delta_{\phi,k}^{\xi}(x,a)$$

(4.9)
$$= M_{\phi,k}\left(\frac{g(x)}{a}\right)$$

$$+ \int_{\mathbb{R}^d} \mathbb{E}\left[\phi\left(\frac{g(x)}{a}|C(\mathbf{0}, \mathcal{H} \cup y)|\right) - \phi\left(\frac{g(x)}{a}|C(\mathbf{0}, \mathcal{H})|\right)\right] dy.$$

We now show that $\Delta_{\phi,k}^{\xi}(x,a)$ reduces to $\Delta_{\phi,k}(g(x)/a)$ defined by (2.8). Put $\beta := g(x)/a$ and $b_d := |B_1^{\mathcal{K}}(\mathbf{0})|$. Since $|C(\mathbf{0},\mathcal{H})| \stackrel{\mathcal{D}}{=} \Gamma_k$, (4.9) yields

$$
\begin{aligned}
&\Delta_{\phi,k}^{\xi}(x,a) - M_{\phi,k}(\beta) \\
&= \int_{\mathcal{K}} \mathbb{E}[(\phi(\beta|C(\mathbf{0},\mathcal{H} \cup y)|) - \phi(\beta|C(\mathbf{0},\mathcal{H})|))\mathbf{1}\{b_d|y|^d \leq \Gamma_k\}]\,dy \\
&= \int_{\mathcal{K}} \mathbb{E}[(\phi(\beta\max(b_d|y|^d, \Gamma_{k-1})) - \phi(\beta\Gamma_k))\mathbf{1}\{b_d|y|^d \leq \Gamma_k\}]\,dy \\
&= \mathbb{E}\int_0^{\Gamma_{k-1}} \phi(\beta\Gamma_{k-1})\,ds + \mathbb{E}\int_{\Gamma_{k-1}}^{\Gamma_k} \phi(\beta s)\,ds - \mathbb{E}[\Gamma_k \phi(\beta\Gamma_k)],
\end{aligned}
$$
(4.10)

where we put $s := |B_{|y|}^{\mathcal{K}}(\mathbf{0})|$. The third term in the right-hand side of (4.10) is

(4.11) $\quad \mathbb{E}[\Gamma_k \phi(\beta\Gamma_k)] = \int_0^{\infty} s\phi(\beta s)\frac{s^{k-1}}{(k-1)!}e^{-s}\,ds = k\mathbb{E}[\phi(\beta\Gamma_{k+1})]$

and likewise, the first term is $(k-1)\mathbb{E}\phi(\beta\Gamma_k)$. Recalling that $\Gamma_k = \sum_{i=1}^{k} \Gamma_{1,i}$ and setting $t = s - \Gamma_{k-1}$, we find that the middle term in the right-hand side of (4.10) is

$$
\mathbb{E}\int_0^{\infty} \phi(\beta(\Gamma_{k-1}+t))\mathbf{1}_{\{t \leq \Gamma_{1,k}\}}\,dt = \mathbb{E}\int_0^{\infty} \phi(\beta(\Gamma_{k-1}+t))e^{-t}\,dt
$$
$$
= \mathbb{E}[\phi(\beta\Gamma_k))] = M_{\phi,k}(\beta).
$$

Combining these expressions for terms in the right-hand side of (4.10) yields

$$
\Delta_{\phi,k}^{\xi}(x,a) = (k+1)M_{\phi,k}(\beta) - kM_{\phi,k+1}(\beta) := \Delta_{\phi,k}(\beta).
$$

The result (2.12) then follows from the proof of Theorem 2.3 of [32]. $\square$

We now extend Theorem 2.1 to cases with $\phi \in \mathcal{F} \setminus \mathcal{F}_0$ (i.e., where $\phi$ "blows up" at 0) via a truncation argument. Given $\varepsilon > 0$, define the functions $\phi^{\varepsilon}:\mathbb{R}^+ \to \mathbb{R}$ and $\phi_{\varepsilon}:\mathbb{R}^+ \to \mathbb{R}$ by

$$
\phi^{\varepsilon}(x) := \begin{cases} \phi(x), & \text{if } x \geq \varepsilon, \\ 0, & \text{otherwise}, \end{cases}
$$
$$
\phi_{\varepsilon}(x) := \begin{cases} \phi(x), & \text{if } x < \varepsilon, \\ 0, & \text{otherwise}. \end{cases}
$$

To prove Theorem 2.1 for $\phi \in \mathcal{F}$ when either $\mathcal{K} = \mathbb{R}^d$ or $d = 1$, we will use the following lemma, whose proof is given in [4].

LEMMA 4.2 ([4], Lemma 5.1). *Given $h \in \mathcal{B}(\mathbb{R}^d)$ and $\delta > 0$, there exists $\varepsilon_0 > 0$ and $n_0 > 0$ such that for $\varepsilon \in (0,\varepsilon_0)$ and $n \geq n_0$ we have $n^{-1}\mathrm{Var}[\langle h, \nu_{n,\phi_{\varepsilon},k}^g\rangle] \leq \delta$.*

Before stating the next lemma, we define for all $\beta > 0$, $y \in \mathbb{R}^d$, $\phi \in \mathcal{F}$ and $\varepsilon > 0$

$$(4.12) \quad \psi(\beta, y) := \mathbb{E}[\phi(\beta|C(\mathbf{0}, \mathcal{H} \cup y)|)\phi(\beta|C(y, \mathcal{H} \cup \mathbf{0})|)] - (\mathbb{E}\phi(\beta\Gamma_k))^2$$

and

$$(4.13) \quad \psi_\varepsilon(\beta, y) := \mathbb{E}[\phi^\varepsilon(\beta|C(\mathbf{0}, \mathcal{H} \cup y)|)\phi^\varepsilon(\beta|C(y, \mathcal{H} \cup \mathbf{0})|)] - (\mathbb{E}\phi^\varepsilon(\beta\Gamma_k))^2.$$

We also define $a_K := \mathbb{E}[\phi^*(\Gamma_k/K)^2] + \mathbb{E}[\phi^*(K\Gamma_k)^2]$ for $K > 0$, and observe for any $K > 1$ that $a_K < \infty$. Also, if $K^{-1} \leq \beta \leq K$, then since $\phi^*$ is decreasing on $(0,1)$ and increasing on $(1, \infty)$,

$$\begin{aligned}
\mathbb{E}[\phi^*(\beta\Gamma_k)^2] &= \mathbb{E}[\phi^*(\beta\Gamma_k)^2 \mathbf{1}\{\Gamma_k \leq 1/\beta\}] + \mathbb{E}[\phi^*(\beta\Gamma_k)^2 \mathbf{1}\{\Gamma_k > 1/\beta\}] \\
&\leq \mathbb{E}[\phi^*(\Gamma_k/K)^2 \mathbf{1}\{\Gamma_k \leq 1/\beta\}] \\
&\quad + \mathbb{E}[\phi^*(K\Gamma_k)^2 \mathbf{1}\{\Gamma_k > 1/\beta\}] \\
&\leq a_K.
\end{aligned}$$
(4.14)

The proof of the next lemma is technical and is given in [4].

LEMMA 4.3. *([4], Lemma 5.2). Let $K > 1$. Then there exists a Lebesgue integrable function $\psi_K^* : \mathbb{R}^d \to [0, \infty)$, such that*

$$(4.15) \quad |\psi_\varepsilon(\beta, y)| \leq \psi_K^*(y) \qquad \forall y \in \mathbb{R}^d \setminus \{\mathbf{0}\}, \varepsilon \in (0,1], \beta \in [1/K, K].$$

Our next two lemmas, proved in detail in [4], show that $V_{\phi^\varepsilon, k}(\beta)$ and $\Delta_{\phi^\varepsilon, k}(\beta)$ defined by (2.9) and (2.8), respectively, converge to $V_{\phi, k}(\beta)$ and $\Delta_{\phi, k}(\beta)$ as $\varepsilon \downarrow 0$.

LEMMA 4.4 ([4], Lemma 5.3). *For all $\beta > 0$ and $k \in \mathbb{N}$, $V_{\phi, k}(\beta)$ satisfies*

$$(4.16) \qquad\qquad \lim_{\varepsilon \downarrow 0} V_{\phi^\varepsilon, k}(\beta) = V_{\phi, k}(\beta).$$

*Moreover, given $K \in [1, \infty)$, it is the case that*

$$(4.17) \qquad \sup\{|V_{\phi^\varepsilon, k}(\beta)| : 0 < \varepsilon \leq 1, 1/K \leq \beta \leq K\} < \infty.$$

LEMMA 4.5. *For any $\beta > 0$ and $k \in \mathbb{N}$, we have*

$$(4.18) \qquad\qquad \lim_{\varepsilon \downarrow 0} \Delta_{\phi^\varepsilon, k}(\beta) = \Delta_{\phi, k}(\beta).$$

*Also, given $K \in [1, \infty)$, it is the case that*

$$(4.19) \qquad \sup\{|\Delta_{\phi^\varepsilon, k}(\beta)| : 0 < \varepsilon \leq 1, 1/K \leq \beta \leq K\} < \infty.$$

PROOF. Let $\phi^*$ be the dominating function given by (2.1). If $\beta > 0$, then it is straightforward (see (5.24) of [4]) to see that $\phi^*(\beta\Gamma_k)^2$ is a nonnegative integrable random variable which dominates $\phi^\varepsilon(\beta\Gamma_k)^2$, so by the dominated convergence theorem, as $\varepsilon \downarrow 0$ we have

$$\mathbb{E}[\phi^\varepsilon(\beta\Gamma_k)] \to \mathbb{E}[\phi(\beta\Gamma_k)]. \tag{4.20}$$

By (2.8) and (4.20), we obtain (4.18). Also, (2.8) implies $|\Delta_{\phi^\varepsilon,k}(\beta)| \leq (k+1)\mathbb{E}[\phi^*(\beta\Gamma_k)] + k\mathbb{E}[\phi^*(\beta\Gamma_{k+1})]$ and the bound (4.19) easily follows from this with (4.14). $\square$

Given $h \in \mathcal{B}(A)$, let $L_h(\phi)$ be the limiting variance in the statement of Theorem 2.1, that is, let

$$L_h(\phi) := \int_A h^2(x) V_{\phi,k}\left(\frac{g(x)}{f(x)}\right) f(x)\, dx$$
$$- \left(\int_A h(x) \Delta_{\phi,k}\left(\frac{g(x)}{f(x)}\right) f(x)\, dx\right)^2. \tag{4.21}$$

LEMMA 4.6. *Given $h \in \mathcal{B}(A)$, it is the case that*

$$\lim_{\varepsilon \downarrow 0} L_h(\phi^\varepsilon) = L_h(\phi). \tag{4.22}$$

PROOF. By assumption, $(g(x)/f(x), x \in A)$ is bounded away from 0 and $\infty$, and $f$ is bounded. Hence by (4.17), the integrand in the first integral in the expression (4.21) for $L_h(\phi^\varepsilon)$ is bounded by a constant, not depending on $\varepsilon$. Similarly, by (4.19), the integrand in the second integral in the expression (4.21) for $L_h(\phi^\varepsilon)$ is bounded by a constant, not depending on $\varepsilon$. By (4.16) and (4.18), for both integrals the integrand converges, as $\varepsilon \downarrow 0$, to the corresponding integrand for $L_h(\phi)$. So, by the dominated convergence theorem, the integrals converge and (4.22) follows. $\square$

PROOF OF THEOREM 2.1. Let $h \in \mathcal{B}(A)$. Given $\delta > 0$, by Lemmas 4.2 and 4.6, we can find $\varepsilon_0 > 0$ and $n_0 > 0$ such that for $\varepsilon < \varepsilon_0$ and $n \geq n_0$, we have

$$|L_h(\phi^\varepsilon) - L_h(\phi)| < \delta \tag{4.23}$$

and $n^{-1}\operatorname{Var}[\langle h, \nu^g_{n,\phi_\varepsilon,k}\rangle] \leq \delta$. The function $\phi^\varepsilon$ lies in the class $\mathcal{F}_0$, so by Proposition 4.2,

$$\lim_{n\to\infty} n^{-1}\operatorname{Var}[\langle h, \nu^g_{n,\phi^\varepsilon,k}\rangle] = L_h(\phi^\varepsilon) \tag{4.24}$$

and hence by the Cauchy–Schwarz inequality, for large enough $n$, we have

$$n^{-1}|\mathrm{Var}\,(\langle h, \nu^g_{n,\phi,k}\rangle) - \mathrm{Var}(\langle h, \nu^g_{n,\phi^\varepsilon,k}\rangle)|$$
$$= |n^{-1}\mathrm{Var}(\langle h, \nu^g_{n,\phi_\varepsilon,k}\rangle) + 2\,\mathrm{Cov}(n^{-1/2}\langle h, \nu^g_{n,\phi^\varepsilon,k}\rangle, n^{-1/2}\langle h, \nu^g_{n,\phi_\varepsilon,k}\rangle)|$$
$$\leq \delta + 2\delta^{1/2}(n^{-1}\mathrm{Var}\langle h, \nu^g_{n,\phi^\varepsilon,k}\rangle)^{1/2} \leq \delta + 2\delta^{1/2}(L_h(\phi) + \delta)^{1/2}.$$

Using (4.23) and (4.24), for large enough $n$, we thus have

$$|n^{-1}\mathrm{Var}[\langle h, \nu^g_{n,\phi,k}\rangle] - L_h(\phi)| \leq 3\delta + 2\delta^{1/2}(L_h(\phi) + \delta)^{1/2}$$

and since $\delta > 0$ is arbitrary, this shows that

$$n^{-1}\mathrm{Var}[\langle h, \nu^g_{n,\phi,k}\rangle] \to L_h(\phi) \qquad \text{as } n \to \infty,$$

which is the first part of the statement of Theorem 2.1.

To prove the stated asymptotic normality of $n^{-1/2}\langle h, \overline{\nu}^g_{n,\phi,k}\rangle$, it suffices to show that for any $h \in \mathcal{B}(A)$,

(4.25) $$n^{-1/2}\langle h, \overline{\nu}^g_{n,\phi,k}\rangle \xrightarrow{\mathcal{D}} N(0, L_h(\phi)).$$

Let $t \in \mathbb{R}$. Set $X_n := n^{-1/2}\langle h, \overline{\nu}^g_{n,\phi,k}\rangle$ and for $\varepsilon > 0$ set $X_n^\varepsilon := n^{-1/2}\langle h, \overline{\nu}^g_{n,\phi^\varepsilon,k}\rangle$. Since $\phi^\varepsilon$ is in $\mathcal{F}_0$, Proposition 4.2 shows that $X_n^\varepsilon \xrightarrow{\mathcal{D}} N(0, L_h(\phi^\varepsilon))$ as $n \to \infty$. Hence,

(4.26) $$\mathbb{E}[\exp(itX_n^\varepsilon)] - \exp(-t^2 L_h(\phi^\varepsilon)/2) \to 0 \qquad \text{as } n \to \infty.$$

Given $\delta > 0$, by Lemmas 4.2 and 4.6, we can choose $\varepsilon > 0$ such that for large $n$,

$$\mathbb{E}[|\exp(itX_n) - \exp(itX_n^\varepsilon)|] \leq \mathbb{E}[|t(X_n - X_n^\varepsilon)|] \leq \delta$$

and also $|e^{-t^2 L_h(\phi)/2} - e^{-t^2 L_h(\phi^\varepsilon)/2}| \leq \delta$ so that combining with (4.26), we have for large $n$ that

$$|\mathbb{E}[\exp(itX_n)] - e^{-t^2 L_h(\phi)/2}| \leq 3\delta$$

and since $\delta$ is arbitrary, this implies (4.25). □

## 5. Proofs of Propositions 2.1 and 2.2.

5.1. *Proof of Proposition 2.1.* First, we identify $V_{\phi,1}(\beta)$ when $\mathcal{K} \neq \mathbb{R}^d$, which implies $-y \notin \mathcal{K}$ for all $y \in \mathcal{K}$. The integral in (2.9) has contributions only from $y \in \mathcal{K}$ and from $\mathbf{0} \in (y + \mathcal{K})$, and these contributions are equal by a symmetry argument. Let $b_d := |B_1^\mathcal{K}(\mathbf{0})|$.

Consider $y \in \mathcal{K}$. Then $|C_1(\mathbf{0}, \mathcal{H} \cup y)|$ has the distribution of $\min(\Gamma_1, b_d|y|^d)$ and $|C_1(y, \mathcal{H} \cup \mathbf{0})|$ has the distribution of $\Gamma_1$, and they are independent. Hence, the integral in (2.9) is equal to

$$2\int_{\mathcal{K}} \mathbb{E}[\phi(\beta\Gamma_1)](\mathbb{E}[\phi(\beta\min(\Gamma_1, b_d|y|^d))] - \phi(\beta\Gamma_1)])\, dy$$

$$= 2M_{\phi,1}(\beta) \int_0^\infty \mathbb{E}[\phi(\beta\min(\Gamma_1, s)) - \phi(\beta\Gamma_1)]\, ds$$

$$= 2M_{\phi,1}(\beta) \mathbb{E} \int_0^{\Gamma_1} (\phi(\beta s) - \phi(\beta\Gamma_1))\, ds.$$

In the last expectation, the first term is equal to $\int_0^\infty \phi(\beta s) P[\Gamma_1 \geq s]$ which comes to $M_{\phi,1}(\beta)$. The second term comes to $M_{\phi,2}(\beta)$ as in (4.11). Thus, the integral in (2.9) is equal to $2M_{\phi,1}(\beta)(M_{\phi,1}(\beta) - M_{\phi,2}(\beta))$ and substituting in (2.9) we find that $V_{\phi,1}(\beta)$ is given by case $k=1$ of formula (2.13) when $\mathcal{K} \neq \mathbb{R}^d$.

To complete the proof of Proposition 2.1, we need to show that $V_{\phi,k}$ is given by (2.13) in the case when $d=1$ (for arbitrary $k$, but still assuming $\mathcal{K} \neq \mathbb{R}^d$). There are only two possibilities for $\mathcal{K}$ and by symmetry it suffices to consider the case with $\mathcal{K} = (0, \infty)$. In this case, the expression (2.9) becomes

$$(5.1) \qquad V_{\phi,k}(\beta) := M_{\phi^2,k}(\beta) + \int_{-\infty}^\infty c^\beta(\mathbf{0}, y)\, dy$$

where

$$c^\beta(\mathbf{0}, y) := \mathbb{E}[\phi(\beta C_\mathbf{0})\phi(\beta C_y)] - (\mathbb{E}[\phi(\beta\Gamma_k)])^2,$$

where $C_\mathbf{0}$ (resp., $C_y$) denotes the length of the $k$-spacing starting at the origin (resp., starting at $y$) with respect to the augmented point set $\mathcal{H} \cup \mathbf{0} \cup y$.

We proceed to evaluate the integral in (5.1). Write $e_k := \mathbb{E}[\phi(\beta\Gamma_k)]$. Then

$$c^\beta(\mathbf{0}, y) = \mathbb{E}[(\phi(\beta C_\mathbf{0})\phi(\beta C_y) - \phi(\beta\Gamma_k)e_k)(\mathbf{1}\{y \leq \Gamma_k\} + \mathbf{1}\{y > \Gamma_k\})]$$

$$= \mathbb{E}[(\phi(\beta C_\mathbf{0})\phi(\beta C_y) - \phi(\beta\Gamma_k)e_k)\mathbf{1}\{y \leq \Gamma_k\}].$$

Integrating over $y$ and setting $\Gamma_0 := 0$, we have that

$$(5.2) \qquad \int_0^\infty c^\beta(\mathbf{0}, y)\, dy = \left(\sum_{j=1}^k I_j\right) - \int_0^\infty \mathbb{E}[\phi(\beta\Gamma_k)e_k\mathbf{1}\{y \leq \Gamma_k\}]\, dy,$$

where we set

$$I_j := \mathbb{E} \int_{\Gamma_{j-1}}^\infty (\phi(\beta C_\mathbf{0})\phi(\beta C_y) \cdot \mathbf{1}\{y \leq \Gamma_j\})\, dy.$$

Recall that $\Gamma_j = \sum_{i=1}^j \Gamma_{1,i}$. We now compute $I_j$ in the case with $1 \leq j \leq k-1$. For such $j$, if $\Gamma_{j-1} < y < \Gamma_j$ then $C_\mathbf{0} = \Gamma_{k-1}$ and $C_y = \Gamma_{j+k-1} - y$;

setting $w = y - \Gamma_{j-1}$, we have for $1 \leq j \leq k-1$ that

$$I_j = \mathbb{E} \int_0^\infty \phi\left(\beta\left(\Gamma_{j-1} + w + (\Gamma_{1,j} - w) + \sum_{i=j+1}^{k-1} \Gamma_{1,i}\right)\right)$$
$$\times \phi\left(\beta\left(\Gamma_{1,j} - w + \sum_{i=j+1}^{j+k-1} \Gamma_{1,i}\right)\right) \mathbf{1}\{\Gamma_{1,j} \geq w\} \, dw.$$

Now take the expectation inside the integral. Since $\Gamma_{1,j}$ is exponential, we have $P[\Gamma_{1,j} \geq w] = e^{-w}$, and by conditioning on this event, using the memoryless property of the exponential distribution and independence of $\Gamma_{1,j}$ from the other random variables in the expression, we obtain

$$I_j = \int_0^\infty \mathbb{E}\phi\left(\beta\left(\Gamma_{j-1} + w + \Gamma_{1,j} + \sum_{i=j+1}^{k-1} \Gamma_{1,i}\right)\right)$$
$$\times \phi\left(\beta\left(\Gamma_{1,j} + \sum_{i=j+1}^{j+k-1} \Gamma_{1,i}\right)\right) e^{-w} \, dw.$$

Now take the integral back inside the expectation. Letting $\Gamma_{1,0}$ be a further independent exponential random variable with density function $e^{-w}, w \geq 0$, we have that

$$I_j = \mathbb{E}\left[\phi\left(\beta\left(\Gamma_{j-1} + \Gamma_{1,0} + \Gamma_{1,j} + \sum_{i=j+1}^{k-1} \Gamma_{1,i}\right)\right)\right.$$
$$\left.\times \phi\left(\beta\left(\Gamma_{1,j} + \sum_{i=j+1}^{j+k-1} \Gamma_{1,i}\right)\right)\right]$$

(5.3)
$$= \mathbb{E}\left[\phi\left(\beta\left(\sum_{i=0}^{k-1} \Gamma_{1,i}\right)\right)\phi\left(\beta\left(\sum_{i=j}^{k+j-1} \Gamma_{1,i}\right)\right)\right]$$
$$= \mathbb{E}[\phi(\beta\Gamma_k)\phi(\beta(\Gamma_{k+j} - \Gamma_j))].$$

To deal with $I_k$, we modify the preceding argument as follows. If $\Gamma_{k-1} < y < \Gamma_k$, then $C_{\mathbf{0}} = y$ and $C_y = \Gamma_{2k+1} - y$. Setting $w = y - \Gamma_{k-1}$, we have that

$$I_k = \mathbb{E} \int_0^\infty \phi(\beta(\Gamma_{k-1} + w))\phi\left(\beta\left(\Gamma_{1,k} - w + \sum_{i=k+1}^{2k-1} \Gamma_{1,i}\right)\right)\mathbf{1}\{\Gamma_{1,k} \geq w\} \, dw.$$

Conditioning on the event that $\Gamma_{1,k} \geq w$ using the memoryless property of the exponential distribution and independence of $\Gamma_{1,k}$ from the other random variables in the expression, we obtain

$$I_k = \mathbb{E} \int_0^\infty \phi(\beta(\Gamma_{k-1} + w))\phi\left(\beta\left(\Gamma_{1,k} + \sum_{i=k+1}^{2k-1} \Gamma_{1,i}\right)\right)e^{-w} \, dw.$$

Letting $\Gamma_{1,0}$ be a further independent exponential random variable, we have that

$$I_k = \mathbb{E}\left[\phi(\beta(\Gamma_{k-1} + \Gamma_{1,0}))\phi\left(\beta\left(\Gamma_{1,k} + \sum_{i=k+1}^{2k-1}\Gamma_{1,i}\right)\right)\right]$$
$$= \mathbb{E}\left[\phi\left(\beta\left(\sum_{i=0}^{k-1}\Gamma_{1,i}\right)\right)\phi\left(\beta\left(\sum_{i=k}^{2k-1}\Gamma_{1,i}\right)\right)\right] = e_k^2.$$

Now as in (4.11) the last term in (5.2) is

$$e_k \int_0^\infty \mathbb{E}[\phi(\beta\Gamma_k)\mathbf{1}\{y < \Gamma_k\}]\,dy = ke_k\mathbb{E}[\phi(\beta\Gamma_{k+1})].$$

Combining this with the preceding expressions for $I_j(j < k)$ and for $I_k$, we may rewrite (5.2) as

$$\int_0^\infty c^\beta(\mathbf{0}, y)\,dy = \left(\sum_{j=1}^{k-1}\mathbb{E}[\phi(\beta\Gamma_k)\phi(\beta\Gamma_{k+j} - \beta\Gamma_j)]\right) + e_k^2 - ke_ke_{k+1}$$
$$= \left(\sum_{j=1}^{k-1}(\mathbb{E}[\phi(\beta\Gamma_k)\phi(\beta\Gamma_{k+j} - \beta\Gamma_j)] - e_k^2)\right) + ke_k(e_k - e_{k+1})$$
$$= ke_k(e_k - e_{k+1}) + \sum_{j=1}^{k-1}\mathrm{Cov}(\phi(\beta\Gamma_k), \phi(\beta\Gamma_{k+j} - \beta\Gamma_j)).$$

By symmetry, for all $\beta$, we have $\int_{-\infty}^0 c^\beta(\mathbf{0}, y)\,dy = \int_0^\infty c^\beta(\mathbf{0}, y)\,dy$, and thus from (5.1) we obtain for all $\beta > 0$ that $V_{\phi,k}(\beta)$ is given by (2.13). This completes the proof of Proposition 2.1.

5.2. *Proof of Proposition 2.2.* We deduce Proposition 2.2 as follows. From the definition (2.9), we obtain

$$V_{\phi,1}(\beta) = M_{\phi^2,1}(\beta) + \int_{\mathbb{R}^d} c(\mathbf{0}, y)\,dy,$$

where $c(\mathbf{0}, y) := \mathbb{E}[\phi(\beta|C_1(\mathbf{0}, \mathcal{H} \cup y)|)\phi(\beta|C_1(y, \mathcal{H} \cup \mathbf{0})|)] - (\mathbb{E}[\phi(\beta\Gamma_1)])^2$. For all $s, t \in \mathbb{R}^+$, let $p(s,t) := P[|C_1(\mathbf{0}, \mathcal{H} \cup y)| > s, |C_1(y, \mathcal{H} \cup \mathbf{0})| > t]$. Then for all $s, t \in [0, |y|^d \omega_d]$, we have

$$p(s,t) = e^{-(s+t)+I(s,t,|y|)}.$$

Otherwise, $p(s,t) = 0$. Hence, for $y \in \mathbb{R}^d$, by the fundamental theorem of calculus, the assumption that $\phi$ is differentiable with $\lim_{t \downarrow 0} \phi(t) = 0$, and

Fubini's theorem,

$$c(\mathbf{0}, y) = \mathbb{E} \int_0^\infty \int_0^\infty \beta^2 \phi'(\beta s)\phi'(\beta t)\mathbf{1}_{\{|C_1(\mathbf{0},\mathcal{H}\cup y)|>s, |C_1(y,\mathcal{H}\cup \mathbf{0})|>t\}}\, ds\, dt$$

$$- \left(\mathbb{E}\int_0^\infty \beta\phi'(\beta u)\mathbf{1}_{\{\Gamma_1>u\}}\, du\right)^2$$

$$= \beta^2 \int_0^\infty \int_0^\infty \phi'(\beta s)\phi'(\beta t)[p(s,t) - e^{-(s+t)}]\, ds\, dt.$$

Since $p(s,t)$ vanishes whenever $(s,t) \notin [0, |y|^d\omega_d]^2$, we obtain

$$c(\mathbf{0}, y) = \beta^2 \int_0^{|y|^d\omega_d}\int_0^{|y|^d\omega_d} \phi'(\beta s)\phi'(\beta t)[e^{-(s+t)+I(s,t,|y|)} - e^{-(s+t)}]\, ds\, dt$$

$$- \beta^2 \int\int_{\max(s,t)\geq |y|^d\omega_d} \phi'(\beta s)\phi'(\beta t)e^{-(s+t)}\, ds\, dt.$$

Therefore,

$$\int_{\mathbb{R}^d} c(\mathbf{0}, y)\, dy$$

$$= \beta^2 \int_{\mathbb{R}^d}\int_0^{|y|^d\omega_d}\int_0^{|y|^d\omega_d} \phi'(\beta s)\phi'(\beta t)[e^{-(s+t)+I(s,t,|y|)} - e^{-(s+t)}]\, ds\, dt\, dy$$

$$- \beta^2 \int_{\mathbb{R}^d}\int\int_{\max(s,t)\geq |y|^d\omega_d} \phi'(\beta s)\phi'(\beta t)e^{-(s+t)}\, ds\, dt\, dy$$

and letting $u := |y|^d\omega_d$, the above becomes

$$\int_0^\infty c(\mathbf{0}, y)\, dy$$

$$= \beta^2 \int_0^\infty \int_0^u \int_0^u \phi'(\beta s)\phi'(\beta t)e^{-(s+t)}[e^{I(s,t,(u/\omega_d)^{1/d})} - 1]\, ds\, dt\, du$$

$$- \beta^2 \int_0^\infty \int\int_{\max(s,t)\geq u} \phi'(\beta s)\phi'(\beta t)e^{-(s+t)}\, ds\, dt\, du.$$

Finally, change the order of integration to obtain

$$\int_0^\infty c(\mathbf{0}, y)\, dy$$

$$= \beta^2 \int_0^\infty \int_0^\infty \phi'(\beta s)\phi'(\beta t)e^{-(s+t)}\int_{\max(s,t)}^\infty [e^{I(s,t,(u/\omega_d)^{1/d})} - 1]\, du\, ds\, dt$$

$$- \beta^2 \int_0^\infty \int_0^\infty \phi'(\beta s)\phi'(\beta t)e^{-(s+t)}\int_0^{\max(s,t)} du\, ds\, dt,$$

which is exactly the desired limit. This proves Proposition 2.2.

**Acknowledgments.** Yu. Baryshnikov thanks the Lehigh Mathematics Department for hosting his research visit. J. Yukich gratefully thanks R. Jiménez for multiple inspiring discussions related to spacings.

Y. BARYSHNIKOV
RM 2C-323
BELL LABORATORIES
LUCENT TECHNOLOGIES
600-700 MOUNTAIN AVENUE
MURRAY HILL, NEW JERSEY 07974
USA
E-MAIL: ymb@research.bell-labs.com

M. D. PENROSE
DEPARTMENT OF MATHEMATICAL SCIENCES
UNIVERSITY OF BATH
BATH BA1 7AY
UNITED KINGDOM
E-MAIL: m.d.penrose@bath.ac.uk

J. E. YUKICH
DEPARTMENT OF MATHEMATICS
LEHIGH UNIVERSITY
BETHLEHEM, PENNSYLVANIA 18015
USA
E-MAIL: joseph.yukich@lehigh.edu